\documentclass[reqno]{amsart}
\usepackage{amssymb}
\usepackage[mathscr]{eucal}

\theoremstyle{plain}
\newtheorem{prop}{Proposition}
\newtheorem{thm}[prop]{Theorem}

\theoremstyle{definition}
\newtheorem*{rem*}{Remark}

\renewcommand{\subset}{\subseteq}

\newcommand{\R}{{\mathbb{R}}}
\newcommand{\x}{{\mathtt{x}}}

\begin{document}

\title
{An observation on positive definite forms}

\author
{Claus Scheiderer}

\subjclass[2010]
{Primary 14\,P\,99;
secondary 11\,E\,25}

\address{Fachbereich Mathematik und Statistik, Universit\"at Konstanz,
  Germany}

\maketitle


\begin{abstract}
Given two positive definite forms $f,\,g\in\R[x_1,\dots,x_n]$, we
prove that $fg^r$ lies in the interior of the sums of squares cone
for large~$r$.
\end{abstract}


\bigbreak
\bigbreak

Fix $n\ge1$ and write $\R[\x]=\R[x_1,\dots,x_n]$. A form (homogeneous
polynomial) $f\in\R[\x]$ is \emph{positive definite} if $f(\xi)>0$
for every $0\ne\xi\in\R^n$. For every integer $d\ge0$ let $\R[\x]_d$
denote the space of forms of degree $d$. Let $\Sigma_{2d}\subset
\R[\x]_{2d}$ be the set of all forms of degree $2d$ that are
\emph{sos}, i.e.\ sums of squares of forms. It is well known that
$\Sigma_{2d}$ is a full-dimensional closed convex cone in
$\R[\x]_{2d}$.

The main result of \cite{sch:ppss} implies (see \cite{sch:ppss}
Remark 4.6):

\begin{thm}\label{ppssthm}%
Let $f,\,g$ be positive definite forms in $\R[x_1,\dots,x_n]$, with
$g$ not constant. Then there is $r_0\ge0$ such that the form $fg^r$
is a sum of squares for all $r\ge r_0$.
\end{thm}

A form $f$ of degree $2d$ will be said to be a \emph{strict sum of
squares}, or \emph{strictly sos}, if $f$ lies in the interior of the
cone $\Sigma_{2d}$. It is equivalent that $f$ has a sum of squares
representation $f=f_1^2+\cdots+f_N^2$ in which $f_1,\dots,f_N$ form a
linear basis of $\R[\x]_d$. The purpose of this note is to observe
that Theorem \ref{ppssthm} can be sharpened as follows:

\begin{prop}\label{result}%
Let $f,\,g$ be positive definite forms in $\R[x_1,\dots,x_n]$, with
$g$ not constant. Then there is $r_0\ge0$ such that the form $fg^r$
is a strict sum of squares for all $r\ge r_0$.
\end{prop}

I am grateful to Amirali Ahmadi for valuable remarks and for
suggesting this question.

Before giving the proof we need a bit of preparation. Given a form
$f\in\Sigma_{2d}$, let $U_f$ be the set of forms $p\in\R[\x]_d$ for
which there exists a real number $c>0$ such that $f-cp^2$ is a sum of
squares. Then $U_f$ is a linear subspace of $\R[\x]_d$, and $f$ is
strictly sos if and only if $U_f=\R[\x]_d$. More generally, the faces
of the cone $\Sigma_{2d}$ are precisely the sets $F_U\Sigma_{2d}:=
\{f\in\Sigma_{2d}\colon U_f\subset U\}$, for $U$ a subspace of
$\R[\x]_d$. (We will not use this fact.) If $U,\,V$ are linear
subspaces of $\R[\x]$, let $UV$ denote the linear subspace spanned by
the products $uv$ with $u\in U$ and $v\in V$. Let $f$ and $g$ be two
forms that are sums of squares. Clearly we have $U_fU_g\subset
U_{fg}$, and $U_f+U_g\subset U_{f+g}$ if $\deg(f)=\deg(g)$. In
particular, the product of two forms that are strictly sos is again
strictly sos.

\begin{proof}
By Theorem \ref{ppssthm} there is $k\ge0$ such that $fg^k$ and
$fg^{k+1}$ are sos. Replacing $f$ with $fg^k$ or $fg^{k+1}$ and $g$
with $g^2$, we can assume that both forms $f$ and $g$ are sos.

Let $\deg(f)=2d$ and $\deg(g)=2e$, write $A=\R[\x]$. Let $q$ be a
strictly sos form with $\deg(q)=2e$ such that $g-q$ is positive
definite, for example  $q=c(x_1^2+\cdots+x_n^2)^e$ for suitable real
$c>0$.
By Theorem \ref{ppssthm} there exists $k\ge1$ such that the forms
$g^k(g-q)$ and $q^k(g-q)$ are both sos. Then for every $r\ge2k$ the
form
$$g^r-q^r\>=\>\sum_{j=0}^{r-1}(g-q)g^jq^{r-1-j}$$
is sos, since this is true for every summand on the right. Since
$q^r$ is strictly  sos, this implies that $g^r$ is strictly sos for
every $r\ge2k$.

Fix $p\in A_d$. There is a real number $c>0$ for which $f-cp^2$ is
positive definite. By Theorem~\ref{ppssthm} there is an integer
$r(p)\ge2k$ such that $(f-cp^2)g^r$ is sos for $r\ge r(p)$. This
implies $U_{p^2}U_{g^r}\subset U_{fg^r}$ for these $r$, and therefore
$pA_{re}\subset U_{fg^r}$ since $g^r$ is strictly sos. Repeat this
argument for every monomial $p$ of degree $d$, and let $r_0$ be the
maximum of the respective numbers $r(p)$. For every $r\ge r_0$ we
then have $A_dA_{re}=A_{d+re}\subset U_{fg^r}$, which means that
$fg^r$ is strictly sos.
\end{proof}

\begin{rem*}
Stengle's Positivstellensatz says that a polynomial $f$ has strictly
positive values if and only if there exist sums of squares of
polynomials $g$ and $h$ such that $fg=1+h$. For a form, being a
strict sum of squares is a certificate for being positive definite.
Therefore Proposition \ref{result} can be seen as a strong
homogeneous version of Stengle's Positivstellensatz, in which the
multiplier can be chosen to be a power of any preassigned
(nonconstant) positive definite form.
\end{rem*}


\end{document}